\documentclass[12pt, reqno]{amsart}
\usepackage{amsmath, amsthm, amscd, amsfonts, amssymb, graphicx, color}
\usepackage[bookmarksnumbered, colorlinks, plainpages]{hyperref}

\theoremstyle{definition}

\theoremstyle{remark}

\numberwithin{equation}{section}

\def\C{{\mathbb C}}

\def\E{{\mathbb E}}

\def\R{{\mathbb R}}

\def\phi{\varphi}

\def\1{{{\ae}}}
\def\2{{{\o}}}
\def\3{{{\aa}}}
\def\6{\, {\rm d}}

\def\<{{\langle}}
\def\>{{\rangle}}

\def\Tr{{\text{\rm Tr}}}

\begin{document}
\setcounter{page}{1}

\title[Uffe Haagerup -- His life and mathematics]{Uffe Haagerup -- His life and mathematics}

\author[M.S. Moslehian, E. St{\o}rmer, S. Thorbj{\o}rnsen, C. Winsl{\o}w]{Mohammad Sal Moslehian$^1$, Erling St{\o}rmer$^2$, Steen Thorbj{\o}rnsen$^3$ $^*$ and Carl Winsl{\o}w$^4$}

\address{$^1$ Department of Pure Mathematics, Center of Excellence in
Analysis on Algebraic Structures (CEAAS), Ferdowsi University of
Mashhad, P.O. Box 1159, Mashhad 91775, Iran.}
\email{\textcolor[rgb]{0.00,0.00,0.84}{moslehian@um.ac.ir}}

\address{$^2$ Department of Mathematics, The Faculty of Mathematics and Natural Sciences, University of Oslo, Norway.}
\email{\textcolor[rgb]{0.00,0.00,0.84}{erlings@math.uio.no}}

\address{$^3$ Department of Mathematics, Ny Munkegade 118, building 1535, 412, 8000 {\AA}rhus C, Denmark.}
\email{\textcolor[rgb]{0.00,0.00,0.84}{steenth@math.au.dk}}

\address{$^4$ Department of Science Education, Faculty of Science, University of Copenhagen, {\O}ster Voldgade 3, 1350 Copenhagen K, Denmark.}
\email{\textcolor[rgb]{0.00,0.00,0.84}{winslow@ind.ku.dk}}

\begin{abstract}
In remembrance of Professor Uffe Valentin Haagerup (1949--2015), as a brilliant mathematician, we review some aspects of his life, and his outstanding mathematical accomplishments.
\end{abstract}



\maketitle

\section{A Biography of Uffe Haagerup}
Uffe Valentin Haagerup was born on 19 December 1949 in Kolding, a
mid-size city in the South-West of Denmark, but grew up in Faaborg
(near Odense). Since his early age he was interested in mathematics. At
the age of 10, Uffe started to help a local surveyor in his work of
measuring land. Soon the work also involved mathematical calculations
with sine and cosine, long before he studied these at school.

\begin{center}
\includegraphics[keepaspectratio=true,scale=0.25]{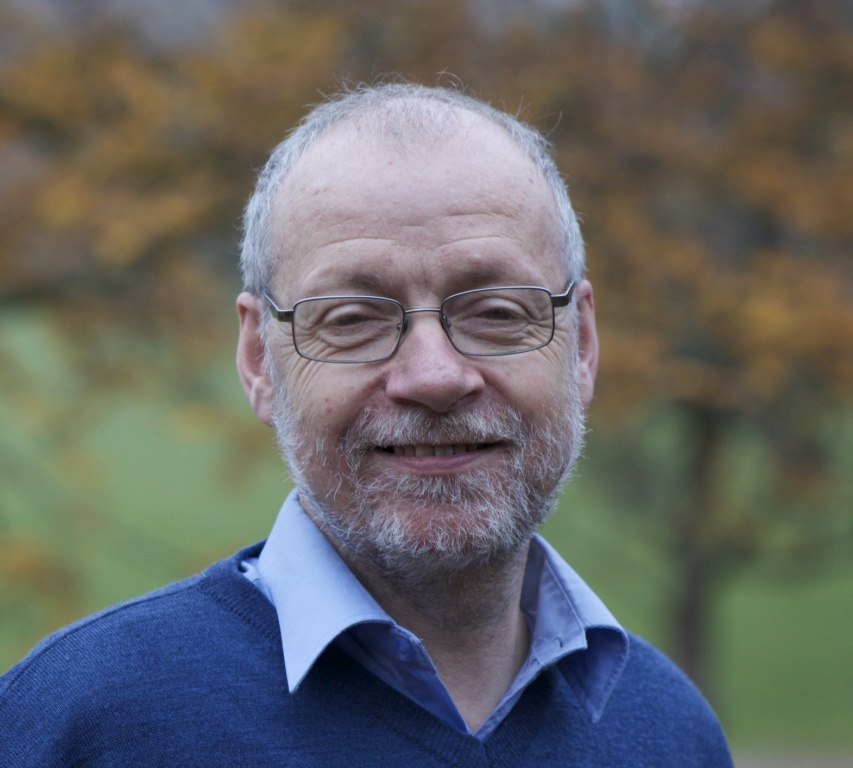}\\
{\textbf{Figure 1.~}Uffe Haagerup -- 2012}
\end{center}

At age 14, Uffe got the opportunity to develop a plan for a new summer
house area close to Faaborg. Due to Uffe's young age, this was
recognized by both local and nationwide media. A plan had previously
been made by a Copenhagen-based engineering company, but their plan
was flawed and eventually had to be discarded.

\begin{center}
\includegraphics[keepaspectratio=true,scale=0.3]{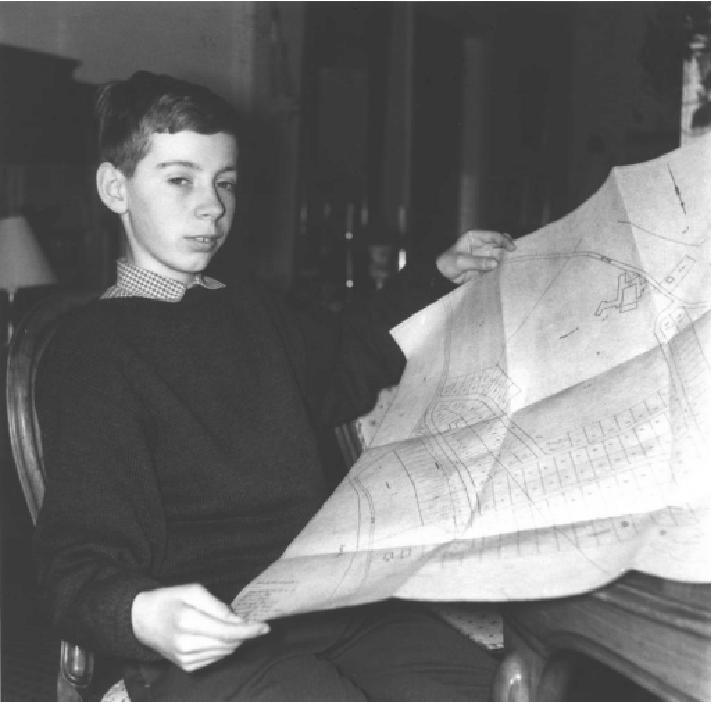}\\
{\textbf{Figure 2.~}A young Uffe Haagerup}
\end{center}
Throughout his childhood, Uffe developed a strong interest in
mathematics, and his skills were several years ahead of those of his
peers. During primary school, he borrowed text books from his
four-year older brother, and thus he came to know mathematics at high
school level. This continued throughout high school, where his
knowledge about mathematics was supplemented by university text
books.

He graduated from high school at Svendborg Gymnasium in 1968. In the
same year, he entered the University of Copenhagen to study
mathematics and physics. He was fascinated by the physical theories of
the 20th century, including Einstein's theory of relativity and
quantum mechanics. His love for the exact language of mathematics led
him to mathematical analysis and in particular the field of operator
algebras, which originally aimed at providing a mathematically exact
formulation of quantum mechanics.

Uffe got his international breakthrough already as a student at the
University of Copenhagen, as he developed an exciting new view on a
mathematical theory developed only a few years before by two Japanese
mathematicians Tomita and Takesaki. From then on, Uffe's name was
acknowledged throughout  the international community of operator
algebraists and beyond.

Uffe received his cand. scient (masters) degree in 1973 from the University of Copenhagen. By the time of his graduation, job opportunities at Danish universities were very limited. Initially Uffe taught a semester at a high school in Copenhagen, but fortunately his talent was recognized by the mathematics department of the newly founded University of Odense (renamed to University of Southern Denmark in 1998), where he was employed from 1974 until his death.

From 1974 to 1977 he served as Adjunkt (Assistant Professor) at the University of Odense. During 1977--79 he had a research fellow position at the same university, which enabled him to devote all of his working hours to scientific work instead of teaching. From 1979 to 1981 he served as Lektor (Associate Professor) and in 1981 the University of Odense promoted him to full Professor at the age of 31, making him the youngest full professor in Denmark.

In 2010-2014, he was on leave from his position at the University of Southern Denmark, to work as a professor at the University of Copenhagen while he held an ERC Advanced grant. In 2015 he returned to his position in Odense. He supervised the following 14 Ph.D. students:\\
\textsl{Marianne Terp (1981), John Kehlet Schou (1991),
Steen Thorbj{\o}rnsen (1998), Flemming Larsen (1999), Jacob v. B. Hjelmborg (2000, co-advisor Mikael R{\o}rdam), Lars Aagaard (2004), Agata Przybyszewska (2006), Hanne Schultz (2006),
Troels Steenstrup (2009), S{\o}ren M{\o}ller (2013), Tim de Laat (2013, co-advisor Magdalena Elena Musat), S{\o}ren Knudby (2014), Kang Li (2015, co-advisor Ryszard Nest),
Kristian Knudsen Olesen (2016, co-advisor Magadalena Elena Musat)}.

His research area mainly falls within operator theory, operator algebras, random matrices, free probability and applications to mathematical physics. Several mathematical concepts
and structures carry his name:\\
\textsc{The Haagerup property} (a second countable locally compact group $G$ is said to have the Haagerup (approximation) property if there is a sequence of normalized continuous
positive-definite functions $\varphi$ which vanish at infinity on $G$ and converge to $1$ uniformly on compact subsets of $G$; see \cite{msm1}), \textsc{the Haagerup subfactor}
and \textsc{the Asaeda--Haagerup subfactor} (Exotic subfactors of finite depth with Jones indices $(5+\sqrt{13})/2$ and $(5+\sqrt{17})/2$;
see \cite{HA}), and \textsc{the Haagerup list} (a list of the only pairs of graphs as candidates for (dual) principal graphs of irreducible subfactors with small index above $4$ and
less than $3+\sqrt{3}$; cf. \cite{msm2}); see also \cite{ALAIN2}.

He spent sabbatical leaves at the Mittag--Leffler Institute in Stockholm, the University of Pennsylvania, the Field Institute for Research in Mathematical Sciences in Toronto and the
Mathematical Science Research Institute at Berkeley.

He served as editor-in-chief for Acta Mathematica from 2000 to
2006. He was one of the editors of the Proceedings of the sixth
international conference on Probability in Banach spaces, Sandbjerg, Denmark, June 16-21, 1986
published by Birkh\"auser in 1990.

He was a member of the ``Royal Danish Academy of Sciences and Letters'' and the ``Norwegian Academy of Sciences and Letters'' and received the following prestigious awards, prizes and honors (\cite{wik, cv}):\\

$\bullet$ The Samuel Friedman Award (UCLA and Copenhagen - 1985) for his solution to the so-called ``Champagne Problem'' posed by Alain Connes.\\

\begin{center}
\includegraphics[keepaspectratio=true,scale=0.22]{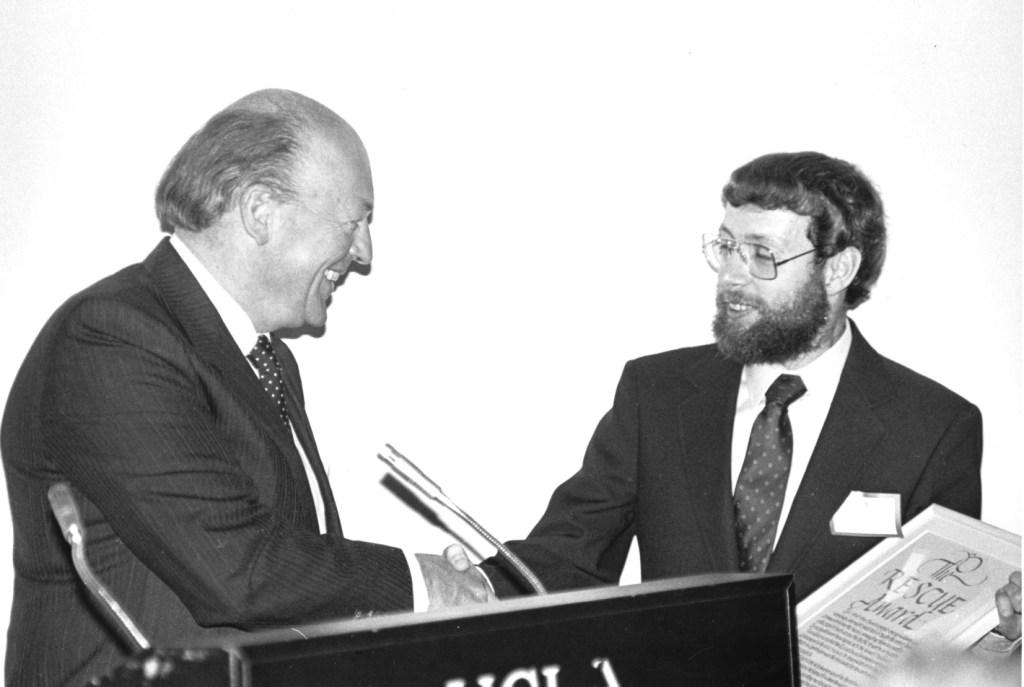}\\
{\textbf{Figure 3.~}Uffe receives the Samuel Friedman Award -- 1985}
\end{center}

$\bullet$ Invited speaker at ICM1986 (Berkeley - 1986).\\

$\bullet$ The Danish Ole R{\o}mer Prize (Copenhagen - 1989). \\

$\bullet$ A plenary speaker at ICM2002 (Beijing - 2002).\\

$\bullet$ Distinguished lecturer at the Fields Institute of Mathematical Research (Toronto - 2007).\\

$\bullet$ The German Humboldt Research Award (M\"unster - 2008).\\

$\bullet$ The European Research Council Advanced Grant (2010--2014).\\

$\bullet$ A plenary speaker at the International Congress on Mathematical Physics ICMP12 (Aalborg - 2012).\\

$\bullet$ The 14th European Latsis Prize from the European Science Foundation (Brussels - 2012) for his ground-breaking and important contributions to operator algebra.\\

$\bullet$ A Honorary Doctorate from East China Normal University (Shanghai - 2013).\\

He also attended numerous conferences and workshops as an invited speaker, such as the 1986 International Congress of Mathematicians in Berkeley, the 2012 International Congress of Mathematicians in Beijing, the 2012 International Congress on Mathematical Physics in Aalborg, and the Conference on Operator Algebras and Applications in Cheongpung.
\begin{center}
\includegraphics[keepaspectratio=true,scale=0.2]{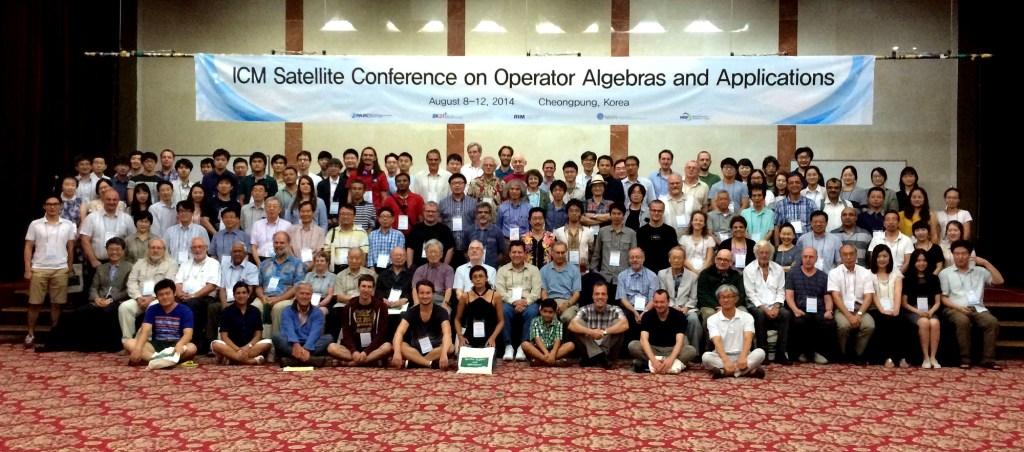}\\
{\textbf{Figure 4.~}Conference in Cheongpung, Korea, 2014}
\end{center}
Uffe had two sons, Peter and S{\o}ren. He tragically drowned on the $5^{{\rm th}}$ of July 2015 while swimming in the sea near Faaborg.\\
\begin{center}
\includegraphics[keepaspectratio=true,scale=0.3]{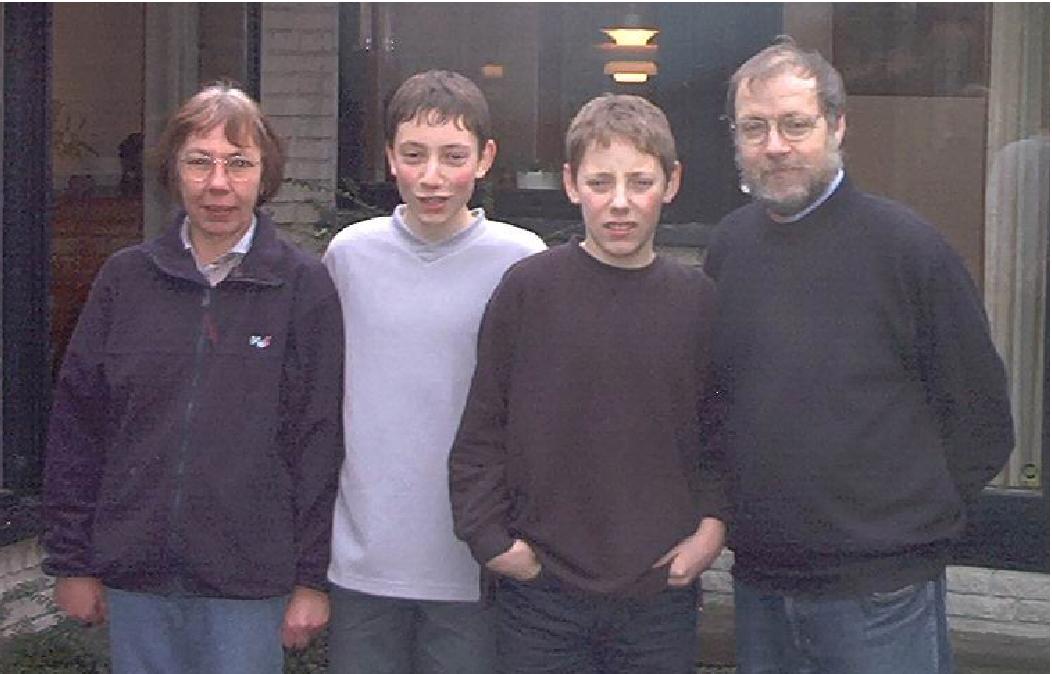}\\
{\textbf{Figure 5.~}Uffe and his family: (From left) Pia, Peter, S{\o}ren, Uffe -- 2002}
\end{center}
In the next two sections, we present some highlights of Uffe Haagerup's mathematical career and works.

\section{Uffe Haagerup's work before 1990}

As mentioned before, Uffe began his studies at the University of
Copenhagen in 1968. At first his main interest was mathematical
physics, and in particular quantum physics.
He got interested in operator algebras via a seminar where papers by
the mathematical physicist Irving Segal, who showed how parts of the
physical theory could be described
by means of operator algebras, were studied.

Let us say a few words on the field operator algebras. This branch of
mathematics was initiated in the years around 1930, when one wished to
develop a mathematical theory
for quantum mechanics, which had been developed a few years earlier by
among others Niels Bohr. The theory was then developed by
mathematicians, in particular John von Neumann, but there were
few active participants in the field until the 1960's. Then physicists
got interested, and the theory of operator algebras became a popular
field.

At the time when Uffe started to learn about operator algebras there was a major breakthrough in the subject. The Japanese mathematician Tomita solved one of the main open problems in
von Neumann algebras, and Takesaki wrote an issue of the Springer Lecture Notes series, which contains the proof plus further developments. Uffe and a fellow student studied these notes in detail.
Then Uffe asked Gert Kj{\ae}rgaard Pedersen if he could write his master thesis on the subject with Pedersen as thesis advisor.
This wish was well received by Pedersen. Uffe wrote his master thesis in the winter 1972-1973, while Pedersen was abroad. It was closely related to Tomita--Takesaki theory and the main results were eventually published in the paper ``The standard form of von Neumann algebras'' \cite{first}, which appeared in 1975.
It is to this day one of his most cited papers and gave him immediately international recognition. His masters thesis also contained another major result, equally published in 1975 \cite{second}: every normal weight on a von Neumann algebra is a supremum of normal states. This solved a problem first formulated by Dixmier.

After this it was unnecessary for Uffe to take a doctoral degree.

In the second half of the 1970's, Uffe produced a number of other important results related to Tomita-Takesaki theory, such as the construction of the $L^{p}$-spaces associated to an arbitrary von Neumann algebra, and a sequence of papers on operator valued weights. One can say that throughout his career, he kept a special affection for (and constantly produced new results in) the area of von Neumann algebras. But he also began, very quickly, to contribute to other fields.

At the University of Odense, Uffe was for a long time the only operator algebraist. However, his colleagues in other fields occasionally told him about famous problems which he was then able to solve. An early example of this was triggered by a problem mentioned to him by his colleague in Banach space theory, Niels J. Nielsen. This gave rise to Uffe's 1978 proof of the best constants in Kinchin's inequality, which consists of 50 pages with difficult classical analysis all the way. Another example is the characterization of simplices of maximal volume in hyperbolic $n$-space, in a joint work with his Odense colleague in algebraic topology, Hans J. Munkholm.

Uffe also contributed to other areas of operator algebra theory, in
particular on $C^*$-algebras, from the late 1970's onwards, and found
new applications in immediately adjacent areas. There is a close
relationship between von Neumann algebras and groups. Groups are
central in mathematics and are algebraic structures where the elements
can be multiplied and have inverses.
 Many constructions of operator algebras involve groups, and the
 algebras often inherit properties
from the underlying group. But the converse is uncommon.
Uffe discovered an example of how properties of groups follow from
operator algebras. He found an example of a so-called non-nuclear
$C^*$-algebra with the metric approximation property.
To do that he started to study a hard analysis problem, and as it
often happened when he solved a problem, he introduced new ideas which
were fruitful for further research.
This time he found a new property of groups, which plays an important
role in geometric group theory. The property is now called the
``Haagerup property'' or a-T-menable, in Gromov's
terminology, as a strong negation of Kazhdan's property (T); see
\cite{COW}.

Uffe didn't forget his background in physics either. A joint work from 1986 with Peter Sigmund, a professor of physics at the University of Odense, shows Uffe's strong analytic powers at work with Bethe's model of energy loss of charged particles as they penetrate matter \cite{sigmund}.
In the fall of 2010, a semester on quantum information theory was held at Institute Mittag--Leffler near Stockholm.
There, Musat came to give a lecture on some joint work with Uffe.
We quote from the report which was written on the program the following year:\\
{\footnotesize``One of the highlights was a pair of visits and talks by Musat, who spoke on her work with
Uffe on factorizable maps, and its implications for the so-called `quantum Birkhoff conjecture', which they showed was false. The first talk generated so much excitement that
questioning went on for more than an hour, with enthusiastic longer
discussions for the rest of Musat's stay."}
\begin{center}
\includegraphics[keepaspectratio=true,scale=0.4]{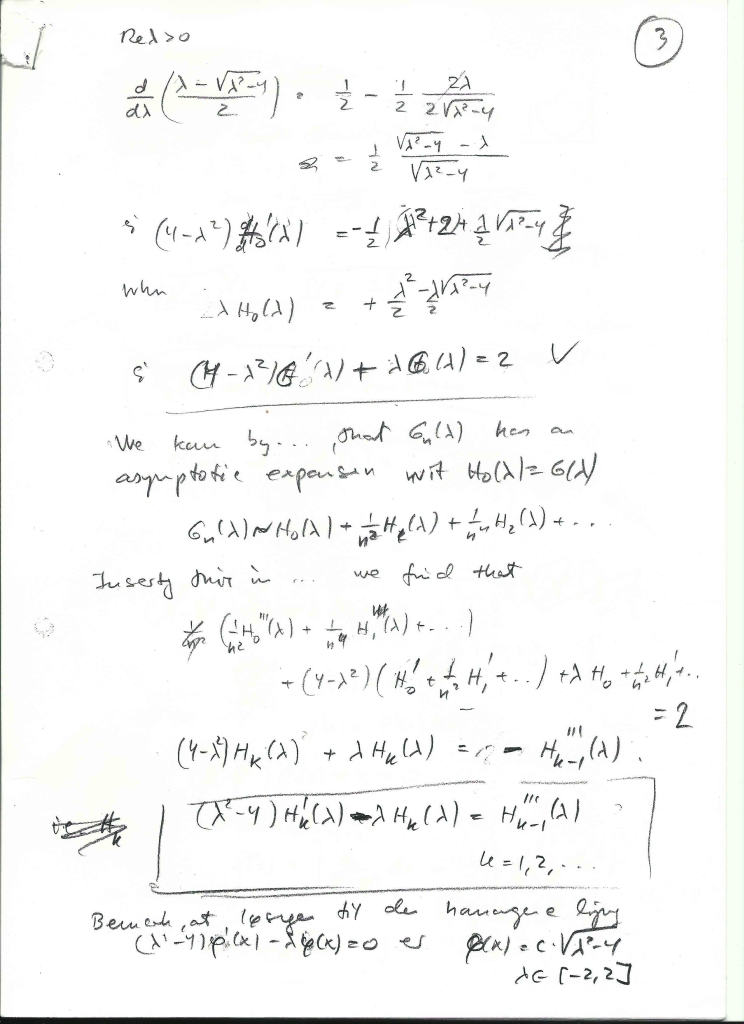}\\
{\textbf{Figure 6.~}A handwritten note by Uffe}
\end{center}
Those of us who have had the pleasure of writing joint papers with Uffe will recognize the following pattern: we had struggled with a problem without success. Then we got into contact with Uffe and told him about the difficulties, where upon he sat down and solved them.

Erling St{\o}rmer recalls one example from a conference in Romania in 1983:\\
{\footnotesize I was going to give a lecture about a formula for the diameter of a set constructed from the states on a
von Neumann algebra. But when I came to the conference I discovered that there were two possible formulas for the diameter, and I was unable to show which one is the correct one.
Fortunately Uffe was there, so I asked him the first day we were there. ``It must be that one", Uffe said and pointed at one of them. In the evening he sat down at his desk,
and the next morning he gave me a 6 page proof showing that the formula he had pointed at, was the right one. Then I could give my lecture with a good feeling.}

A couple of years later, in 1985,  St{\o}rmer and Haagerup shared an apartment in Berkeley in California, for a month. They followed up their work with the diameter formula and ended up with an 80 pages long
paper. In addition to learning much mathematics from this collaboration, St{\o}rmer learned one more thing, namely patience. He wrote a draft, which he sent to Uffe early in the fall of 1987.
But Uffe got ill that fall, so he was delayed in the work of finishing the manuscript, but it took a long time for other reasons too, because Uffe was a very patient mathematician, who could
keep a manuscript in his drawer for a long time before he had them typed and published. Some he never published but sent copies of them to his colleagues. So it was far into the winter before he
gave the final manuscript to the secretary who was going to type it for them. But that also took a long time, so it was only sent to a journal late in the following summer, and then it took at least another two years before it was finally published.

Uffe spent the academic years 1982-83 in the USA, first at UCLA, then at the University of Pennsylvania, both places simultaneously with the 3 years younger Vaughan Jones. At that time Jones
showed some very important results on von Neumann algebras. They were about factors, which in a way are the building blocks in the theory, and have the property in common with the $n\times n$
matrices that their center, i.e. the operators in the algebra which commutes with all the others, consists only of the scalar multiples of the identity operator. For an inclusion
$A\supseteq B$ of factors, Jones introduced an index, which in a way measures the difference of the sizes of the two factors. For some factors he found a formula for the index,
which turned out to be very important for the theory of knots. This was a sensation, as it was an application of the infinite dimensional theory of von Neumann algebras to the finite
dimensional knot theory.
At the world congress in mathematics in 1990 Jones was rewarded the Field's Medal, which is the most prestigious award a mathematician younger than 40 years can get. We return to some of Uffe's contributions to subfactor theory in the next section.

We have now arrived at Uffe's most famous result. He himself also considered this to be his best result ever. When we indicated what von Neumann algebras are, we started with the $n\times n$
matrices. Consider an infinite long increasing sequence of matrix algebras, where each matrix algebra contains the previous ones. From this infinite sequence one can generate many different
von Neumann algebras, and in particular factors by use of states. They are called hyperfinite or injective factors, and have been central in the theory since Murray and von Neumann started
the development of the theory in the 1930's. Factors are divided into classes of types ${\rm I}, {\rm II}$ and ${\rm III}$, and each of these has several subclasses. In this connection type III is most important, and this class is further divided into the types ${\rm III}_{\lambda}$, where $\lambda$ moves through the interval from $0$ to $1$. These are the most ``infinite" von Neumann algebras, and were considered almost untractable before Connes' groundbreaking work, based on Tomita-Takesaki theory.

Alain Connes received the Field's Medal in 1982 for his seminal work
on von Neumann algebras, and especially for his classification of
hyperfinite factors, published in 1976. In particular, he obtained a
complete classification of hyperfinite factors of type  ${\rm
  III}_{\lambda}$, where $0\leq\lambda < 1$.
But there was one problem he did not succeed to solve, namely whether
there is one or more hyperfinite type ${\rm III}_1$-factors. Uffe
visited Connes in 1978 at his country house in Normandie,
where they discussed the problem. Hjelmborg, while preparing an
interview with Uffe in 2002, got the following description from Connes
on these discussions:\\
{\footnotesize ``We had long and intense discussions in my country
  house ending up when both of us got a terrible migraine
  \cite{interview}''.}
\begin{center}
\includegraphics[keepaspectratio=true,scale=0.3]{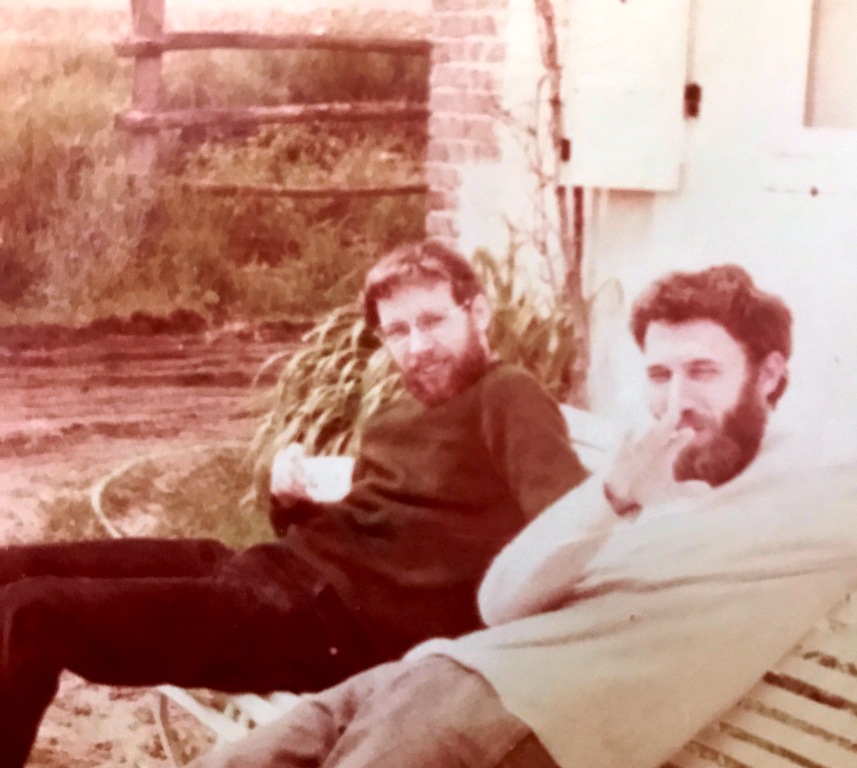}\\
{\textbf{Figure 7.~}Uffe Haagerup (left) and Alain Connes (right), about 40 years ago}
\end{center}
Uffe thought much about the problem later on, but didn't get the opportunity to work seriously on it before the years 1982-83. Based on Connes' work, Uffe finally solved the problem in the fall of
1984, by showing that there is only one hyperfinite type ${\rm
  III}_1$-factor. The proof was published in Acta Mathematica in 1987
and was over 50 pages long (see \cite{haa2}). It demonstrated convincingly
how exceptionally good Uffe was in analysis.
The problem was known as the ``Champagne Problem'', as Connes had promised a fine bottle of Champagne to the person who could solve it. Uffe received the announced Champagne from Connes for the result, as well as the Samuel Friedman Award in 1985.
In an obituary written shortly after Uffe's passing \cite{ALAIN}, Connes expressed his admiration as follows: \\
{\footnotesize
Uffe Haagerup was a wonderful man, with a perfect kindness and openness of mind, and a mathematician of incredible power and insight.
His whole career is a succession of amazing achievements and of decisive and extremely influential contributions to the field of operator algebras, $C^*$-algebras and von Neumann algebras.
(...)
From a certain perspective, an analyst is characterized by the ability of having ``direct access to the infinite'' and
Uffe Haagerup possessed that quality to perfection. His disappearance is a great loss for all of us.}
\begin{center}
\includegraphics[keepaspectratio=true,scale=0.35]{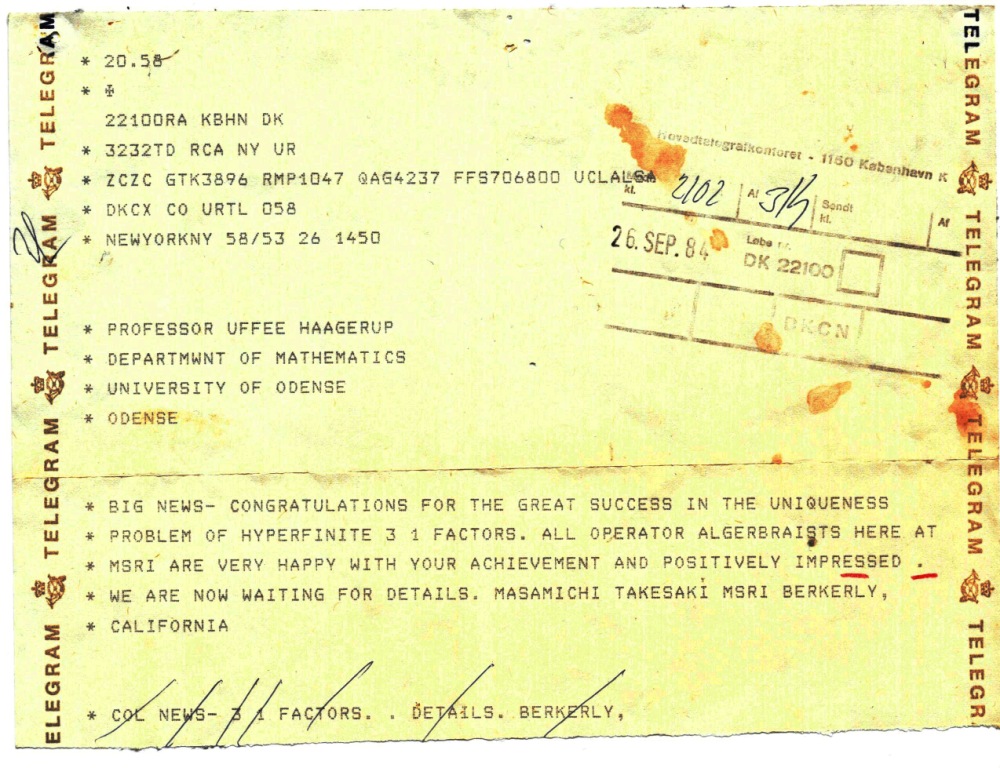}\\
{\textbf{Figure 8.~}Congratulation telegram from Masamichi Takesaki, for the solution of the Champagne Problem}
\end{center}

\begin{center}
\includegraphics[keepaspectratio=true,scale=0.3]{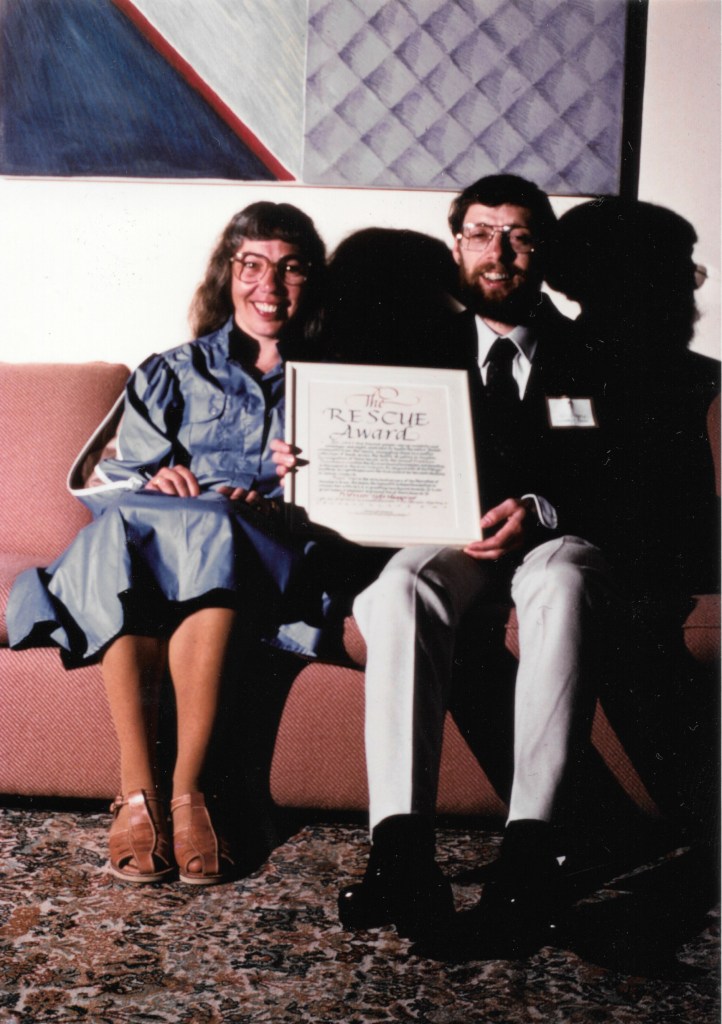}\\
{\textbf{Figure 9.~}Uffe and his wife Pia in 1985}
\end{center}

\section{Uffe Haagerup's work after 1990}

As appears from the preceding section, most of Uffe's research was
centered around the theory of von Neumanns algebras, which he mastered
to the highest international level, and in particular he
 became known for using von Neumann algebra techniques in order to
 solve $C^*$-algebra problems and more generally for using methods
 from analysis to prove results that had
been established previously by other methods. A couple of examples of
this are the following:

\begin{itemize}

\item[(a)] In the paper
``Random Matrices with Complex Gaussian Entries''(ref.~\cite{ht1})
new proofs were
given for the limiting behavior of the empirical spectral
distribution and the smallest and largest eigenvalues of certain
Gaussian random matrix ensembles. In particular these results include
the celebrated semi-circle law of Wigner (see \cite{Wi}). Where
previous proofs of the mentioned
results involved a substantial amount of combinatorial work, Uffe took
the point of view of studying the ``moment generating function''
$s\mapsto \E[\Tr(\exp(sA))]$, where $A$ is the random matrix under
consideration, $\E$ denotes expectation and $\Tr$ denotes the trace.
Expanding this function as a power series, Uffe and his co-author could
identify it explicitly in terms of certain
hypergeometric functions. This approach resembles methods from
analytic number theory, which Uffe was actually quite interested in
and taught several courses on.

\item[(b)]
Another example is the paper ``On Voiculescu's $R$- and $S$-transforms
for free non-commuting random variables'' \cite{haa3} in which (among
other results) Uffe provided a new and completely
analytical proof of the additivity
(with respect to free convolution) of Voiculescu's $R$-transform (see
\cite{Vo}). Voiculescu's original proof was based on the
Helton--Howe formula from representation theory, and other proofs
(e.g.\ by Nica and Speicher; see \cite{NiSp}) are based on the development
of some rather heavy combinatorial machinery.
Uffe's proof is based on Banach-algebra techniques, which
he used e.g.\ to express the $R$-transform explicitly as an analytic
function in a neighborhood of zero. Voiculescu
recently used Uffe's approach to establish a key formula for the
analog of the $R$-transform in Voiculescu's recent theory of bi-free
probability. As it happens, neither Voiculescu's original approach,
nor the combinatorial approach work in the bi-free setting.
In his talk at the celebration of Uffe's 60'th birthday,
Voiculescu gave the following general characterization of Uffe's papers
(quoted freely from memory): ``Everything is very clear and looks very
easy. Then suddenly a `miracle' occurs, from which everything falls out and
is again clear and easy''.
\end{itemize}

In the rest of this section we outline a few highlights from the second
half of Uffe's career. They are listed in chronological
order and should be viewed mainly as examples of his impressive
achievements. Many other of Uffe's results from the last 25 years
equally deserve to be highlighted, but time and space limitations
prevents a thorough encyclopeadic approach.

\begin{itemize}

\item Uffe was inspired by the time he spent with with Jones in 1983 (cf. previous section) and started to work on subfactor theory, as it is called. He eventually ended up by solving a central problem in the
theory, which Jones had left open, namely to find a finite depth,
irreducible subfactor of the hyperfinite factor of type II$_{1}$ with
index strictly between 4 and $3+\sqrt{2}$. Haagerup proved that no
such subfactor can have index smaller than $(5+\sqrt{13})/2$
(\cite{msm2}), and subsequently, with Asaeda,  (\cite{HA}), he proved
the existence and uniqueness of a (finite depth, irreducible)
subfactor of precisely this index. This subfactor, called the
``Asaeda--Haagerup subfactor'', has a very complicated construction
and cannot be constructed by standard methods. Although this result is
less spectacular than the uniqueness of the hyperfinite type
III$_{1}$-factor, it is certainly a second major problem left open by
a Field-medalist and solved by Uffe.

\item
From around the late 1990's Uffe (and collaborators)
made important
contributions to Voiculescu's free probability theory. In
ref.~\cite{ht5} he proved (jointly with Thorbj{\o}rnsen) that
the operator norm of a non-commutative polynomial in several
independent GUE-random matrices converges
almost surely, as the dimension goes to infinity, to the
limit anticipated by free probability theory. This further lead to the
settlement (in the positive) of the conjecture on the existence of
non-invertible
elements in the extension semi-groups of the reduced $C^*$-algebras
associated to the free groups. Jointly with Schultz, he also made
huge progress on the invariant subspace problem. Specifically they proved
in \cite{HS}
that any operator $T$ in a II$_1$-factor has a non-trivial invariant
subspace affiliated with the von Neumann algebra generated by $T$,
provided that the Brown measure of $T$ is non-trivial.

\item In 2008 Uffe and Musat solved in \cite{HaMu} a long standing
 conjecture by Effros--Ruan and Blecher by establishing the following
 Grothendieck type inequality: For any $C^*$-algebras $A$ and $B$ and
 any jointly completely bounded bilinear form $u\colon A\times B\to\C$ there exist states $f_1,f_2$ on $A$ and $g_1,g_2$ on $B$,
 such that
\[
|u(a,b)|\le\|u\|_{\rm jcb}\big(f_1(aa^*)^{1/2}g_1(b^*b)^{1/2}
+f_2(a^*a)^{1/2}g_2(bb^*)^{1/2}\big),
\]
for any $a$ in $A$ and $b$ in $B$. The \emph{jointly completely bounded norm}
$\|u\|_{\rm jcb}$ may be defined as the completely bounded norm of the
mapping $A\to B^*$ associated to $u$. The work of Uffe and Musat
extended previous work by  Pisier and  Shlyakhtenko (see \cite{Pi_Sh}).

\item In a series of two papers (\cite{H_dL1},\cite{H_dL2})
Uffe and de Laat showed recently that all connected, simple Lie groups with
real rank greater than or equal to
 2 do not have the Approximation Property (AP) (see e.g.\
 \cite{H_dL1} for the definition of this property).
 Since connected, simple
 Lie groups with real rank 0 (resp.\ 1) are known to be amenable
 (resp.\ weakly amenable), and since amenability implies weak
 amenability, which again implies (AP), Uffe and de Laat's result
 shows that connected simple Lie groups have (AP), if and only if their
 real rank is at most 1.
 Specifically Uffe and de Laat proved that the symplectic group
$\textrm{Sp}(2,\R)$ and its universal covering group
$\widetilde{\textrm{Sp}}(2,\R)$ do not have the (AP). A few years
before it had been established by Lafforgue and de~la~Salle that
$\textrm{SL}(3,\R)$ does not have the approximation property (see
\cite{L_dlS}). Furthermore it is well-known that any
connected simple Lie Group with real rank
greater than or equal to 2 has a closed connected subgroup, which is
\emph{locally} isomorphic to either or $\textrm{Sp}(2,\R)$ or
$\textrm{SL}(3,\R)$, and hence isomorphic to a quotient of one of the
universal covering groups $\widetilde{\textrm{Sp}}(2,\R)$
or $\widetilde{\textrm{SL}}(3,\R)$ by a discrete normal central subgroup.
Combining the results mentioned above, Uffe and de~Laat's result may
then be deduced from the fact that (AP) passes from a group to its
closed subgroups.

\item
In recent years Uffe became interested in the famous problem on the
possible amenability of the smallest of the Thompson groups, here denoted by
$F$. In 2015 he published the joint paper ref.~\cite{HHR} with
Ramirez--Solano and his youngest son,
S{\o}ren, in which they give precise lower bounds for the norms of two
operators associated to the generators of $F$. By work
of Kesten, the amenability of $F$ is equivalent to the
statement that these norms equal 3 and 4, respectively. Extensive computer
calculations, performed by Uffe and his co-authors, suggest that the norms
are approximately around 2.95 and 3.87, respectively, but their upper
bounds are not precise enough to establish non-amenability. In the
paper \cite{H_KO}, Uffe and Knudsen Olesen established
that if the reduced $C^*$-algebra of the larger Thompson group, $T$, is
simple, then $F$ is non-amenable. Very recently Le Boudec and Matte
Bon proved that non-amenability of $F$ is in fact equivalent to
simplicity of $C^*_{\textrm{r}}(T)$ (see \cite{LB_MB}).

\end{itemize}

\section{Uffe Haagerup as teacher and supervisor}

Many of the numerous students who were taught by Uffe over the years at the University of
 Southern Denmark, mainly saw him as someone who was able to write
 incredibly fast (while still producing readable text) on a
 blackboard. Little did they realize that they were enjoying the
 privilege of being lectured to by one of the greatest and most
 influential Danish mathematicians of all times. Their ignorance is (partly)
 excused by Uffe's general attitude and appearance, to which the
 word ``modest'' immediately springs to the mind of anyone who have
 met him. Of course the students who took more advanced courses with Uffe,
 and in particular those who wrote their masters or Ph.D-thesis under
 his supervision, eventually realized that there was full concordance
 between the pace of his handwriting and that of his mathematical
 mind. One of Uffe's students (Carl Winsl{\o}w) at the University of Southern Denmark
remembers Uffe's marvelous teaching and supervision as
follows:

{\footnotesize My first memories of Uffe Haagerup date back to a linear algebra
course in the late 1980's, at the University of Odense. The lectures
were astonishing, superior to all other I have attended. While his
teaching was spontaneous (no manuscript) and very lively, leaving the
audience in no doubt on the rationale for the current details, he
filled the blackboards with crystal clear proofs and simple examples
-- always more elegant and illuminating than those in the
textbook we had. He repeated the same act in later courses I had the
chance to take with him, on functional analysis, von Neumann algebras
and so on.

Later, at weekly meetings with him as my master thesis
supervisor, the blackboard was replaced with his favorite working
instrument: blanksheets of paper and a classical pencil, which was frequently
sharpened while the sheets where filled, and the sheets were
eventually stapled when some proof was done. My thesis was to be an
exposition of the details of Connes' 1973 paper \cite{ALAIN3}. Of
course Uffe knew this monumental
work intimately; in fact one of his most famous achievements was to
complete the classification in question by proving the uniqueness of
the injective type III$_{1}$ factor, in 1984. At the supervision
meetings, the following often happened: I had struggled with some
elegant but very short proof from Connes' paper, and asked Uffe about
it. He would take a look at the French text, mainly to get the result
to be proved, then provide an elaborate and crystal clear proof on
white sheets, out of his head, which I suspect was quite independent
from the explanation in the paper. It also happened, sometimes, that I
brought up other questions for which I could not find an answer in the
literature. Usually, he would go: ``Yes, I once thought about that'',
take a stack of white paper, and begin writing a sequence of lemmas
and so on - often quite technical with subtle inequalities that were
stated without hesitation and then proved quickly, with occasional
corrections done by simply barring a line or too (I don't recall him
having an eraser). On seldom, happy occasions, he would reach for one
of his endless folders of stapled, handwritten manuscripts, which filled
the shelves in his office - but even then, he usually ended up writing
a new one from scratch.}

This little anecdote is communicated here because we think any of his
students (graduate or undergraduate) would recognize the point: Uffe
incarnated mathematical creativity in a way that is shared by few (if
any) they have met. For him there was a perfect continuity between
``teaching'' and ``research'' - it was about producing and sharing
mathematical ideas. Even in his lectures on linear algebra (where, of
course, no results were new) one got an experience how reasoning and
connections are built ``in vivo''.

In the literature on the modes and effects of mathematics teaching,
the activities in which mathematicians build new knowledge have
sometimes been used as an ideal model for the activity of the student;
the teacher, then, should arrange situations in which the student
could learn by solving and posing problems. Uffe certainly practiced
this art in many ways. But his acts of ``direct teaching'' (allured to
above) were also very far from the caricature image that is sometimes
presented as the ``opposite'' of that ideal: lectures which leave the
students completely passive. Indeed, many lectures fail to help
students to go beyond the role of spectator. But, as his students would say, Uffe's did not.

In the first decades of his academic career,
Uffe only took on a single Ph.D-student: Marianne Terp. From around the mid
 1990's he changed his policies on this matter, partly influenced by
 general tendencies at the Danish Universities, and until his
 death he acted as supervisor on at least another 13
 Ph.D-theses. He never obtained a Ph.D-degree himself; a fact that was used as a
 friendly (and absurd) tease among students and colleagues. In 2013 he
 could, however, put an end to the teasing, as he was awarded an
 honorary doctoral degree from East China Normal University.
\begin{center}
\includegraphics[keepaspectratio=true,scale=0.80]{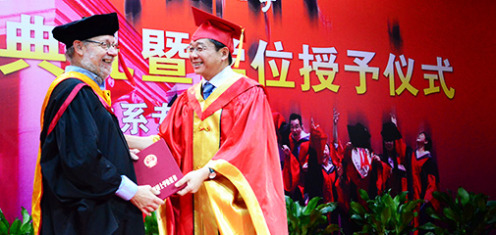}\\
{\textbf{Figure 10.~}Uffe Haagerup and ECNU President Qun Chen at the award ceremony -- 2013}
\end{center}
 In the minds of all his students and post docs, Uffe will always
 stand out as a true master of mathematical thinking and a great
 source of inspiration.
 Collaborating with Uffe was an immense privilege, and his modest and
 kind personality neutralized the feeling of mathematical inferiority
 one could easily get stung by in his presence. Arrogance was simply
 not a part of his character. A very precise description of Uffe as
 teacher and supervisor can be expressed with the
 Japanese term, {\em
 sensei}. It can be used to translate a variety of English terms:
teacher, master, professor, expert, senior. Literally, it means ``the one who
proceeds'' (or walks ahead of) you.
\begin{center}
\includegraphics[keepaspectratio=true,scale=0.2]{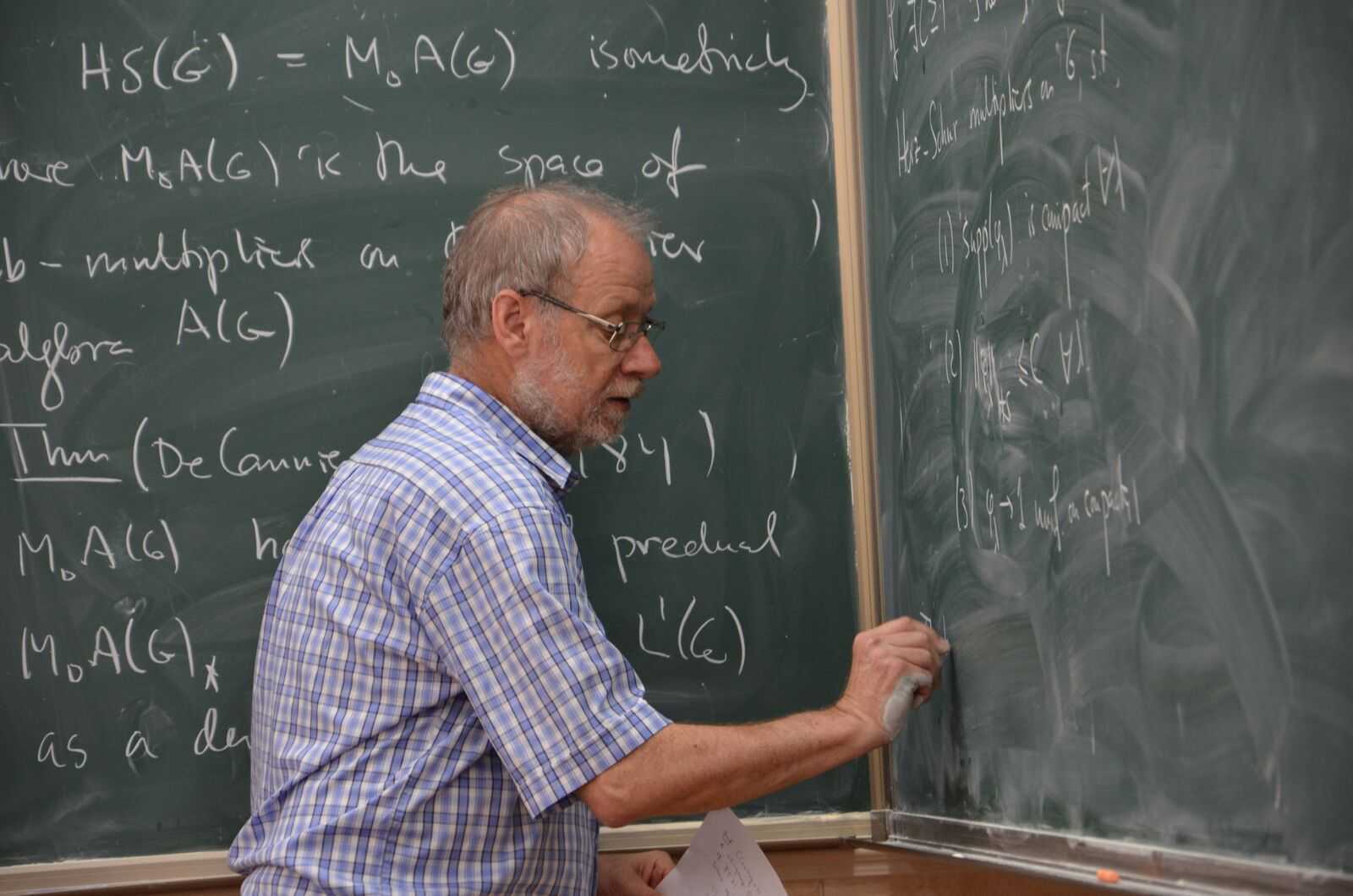}\\
{\textbf{Figure 11.~}Uffe Haagerup in action}
\end{center}

Uffe was a sensei in all the meanings of the word. A sensei badly missed, but whose memory is gladly honored.

\section{Bibliometrics}

Utilizing MathSciNet (MR)\cite{mr}, Zentralblatt MATH (Zbl) \cite{zbl}, Scopus \cite{scopus} and Web of Science (WOS) \cite{wos}, we present some quantitative analysis of Uffe's publications until.

The first three most cited publications of Uffe in MR are:

$\bullet$ U. Haagerup, An example of a nonnuclear $C^*$-algebra, which has the metric approximation property. Invent. Math. 50 (1978/79), no. 3, 279-293. (206 citations)\\

$\bullet$ Michael Cowling and U. Haagerup, Completely bounded multipliers of the Fourier algebra of a simple Lie group of real rank one. Invent. Math. 96 (1989), no. 3, 507-549. (114 citations)\\

$\bullet$ Jean De Canni\'ere and U. Haagerup, Multipliers of the Fourier algebras of some simple Lie groups and their discrete subgroups. Amer. J. Math. 107 (1985), no. 2, 455-500. (111 citations)\\

The first three most cited publications of Uffe in ZbMath are:\\

$\bullet$ U. Haagerup, An example of a non nuclear $C^*$-algebra, which has the metric approximation property, Invent. Math. 50, 279-293 (1979). Zbl 0408.46046 (138 citations)\\

$\bullet$ U. Haagerup, The standard form of von Neumann algebras, Math. Scandinav. 37(1975), 271-283 (1976). Zbl 0304.46044 (102 citations)\\

$\bullet$ Michael Cowling and U. Haagerup, Completely bounded multipliers of the Fourier algebra of a simple Lie group of real rank one, Invent. Math. 96, No.3, 507-549 (1989). Zbl 0681.43012 (81 citations)\\

The first three most cited publications of Uffe in WOS are:\\

$\bullet$ U. Haagerup, An example of a nonnuclear $C^*$-algebra, which has the metric approximation property. Invent. Math. 50 (1978/79), no. 3, 279-293. (297 citations)\\

$\bullet$ U. Haagerup, The standard form of von Neumann algebras, Math. Scandinav. 37(1975), 271-283. (232 citations)\\

$\bullet$ Michael Cowling and U. Haagerup, Completely bounded multipliers of the Fourier algebra of a simple Lie group of real rank one. Invent. Math. 96 (1989), no. 3, 507-549. (143 citations)\\

The first three most cited publications of Uffe in Scopus are:\\

$\bullet$ U. Haagerup, An example of a nonnuclear $C^*$-algebra, which has the metric approximation property. Invent. Math. 50 (1978/79), no. 3, 279-293. (268 citations)\\

$\bullet$ Michael Cowling and U. Haagerup, Completely bounded multipliers of the Fourier algebra of a simple Lie group of real rank one. Invent. Math. 96 (1989), no. 3, 507-549. (136 citations)\\

$\bullet$ U. Haagerup, All nuclear $C\sp{\ast} $-algebras are amenable. Invent. Math. 74 (1983), no. 2, 305--319. (118 citations)\\

The number of Uffe's publications recorded in MR and Zbl are 106 and 109, respectively. According to MR, they are cited 2407 times by 1250 authors. Functional analysis is the subject where Uffe has published most of his articles and where there are most citations to Uffe's works.

According to Zbl, the first three journals with most of Uffe's publications are Journal of Functional Analysis (13 papers), Duke Mathematical Journal (6 papers) and Mathematica Scandinavica (6 papers). He has had 52 collaborators; among them E. St{\o}rmer, S. Thorbj{\o}rnsen and K. J. Dykema with 9, 7 and 5 papers, respectively, have most joint papers with him.

Web of Science records 86 publications by Uffe. The sum of the times his papers are cited is 2775 and without self-citations is 2631. The average citation per publication is 32.27, and Uffe's h-index is 28.\\
\begin{center}
\includegraphics[keepaspectratio=true,scale=0.4]{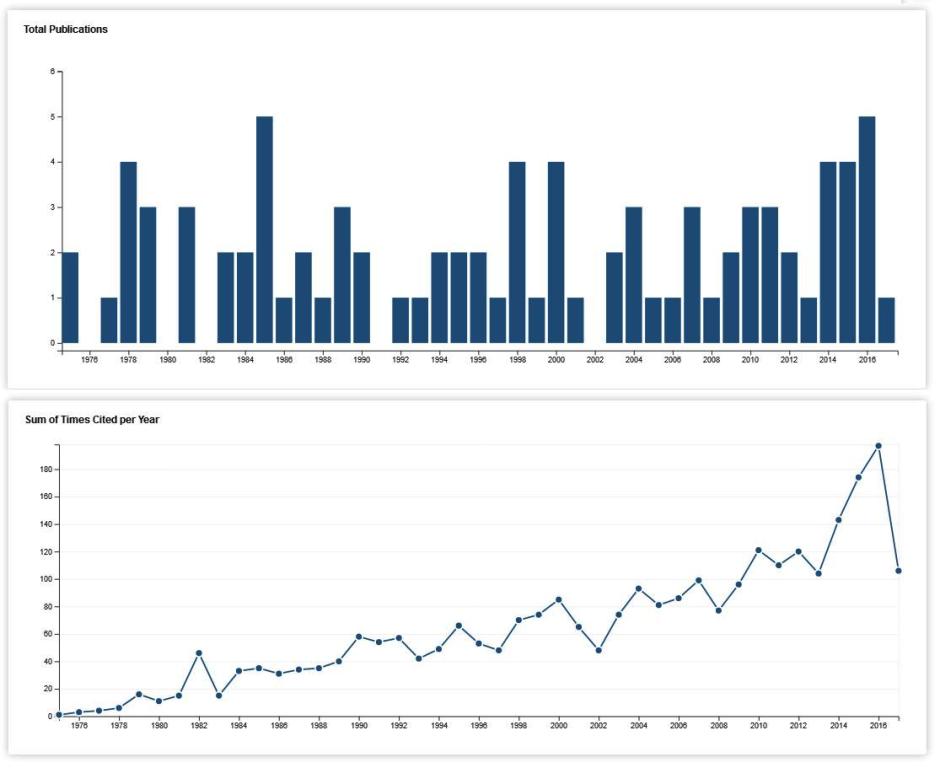}\\
{\textbf{Figure 12.~}Bibliometrics - Ref. Web of Science - Sep 8, 2017}
\end{center}

Scopus presents 69 document for Uffe Haagerup. Records show 1877 total citations by 1431 documents for him. His scopus h-index is 22.

His first paper appearing in MathSciNet is

$\bullet$ Haagerup, Uffe Normal weights on $W^{\ast}$-algebras. J. Functional Analysis 19 (1975), 302--317.

and his last sole author paper is

$\bullet$ Haagerup, Uffe On the uniqueness of the injective $\textrm{III}_1$ factor. Doc. Math. 21 (2016), 1193--1226.

which is typed by Hiroshi Ando and completed (due to some missed pages of the original handwritten note) by Cyril Houdayer and Reiji Tomatsu after Uffe passed away.

\section{Publications by Uffe Haagerup}
His papers listed in MathSciNet are as follows:\\

\begin{itemize}
\item Haagerup, Uffe; Olesen, Kristian Knudsen Non-inner amenability of the Thompson groups T and V. J. Funct. Anal. 272 (2017), no. 11, 4838-4852.

\item Haagerup, Uffe On the uniqueness of the injective $\textrm{III}_1$ factor. Doc. Math. 21 (2016), 1193-1226.

\item Haagerup, Uffe; Knudby, S{\o}ren; de Laat, Tim A complete characterization of connected Lie groups with the approximation property. Ann. Sci. Ec. Norm. Super. (4) 49 (2016), no. 4, 927-946.

\item Ando, Hiroshi; Haagerup, Uffe; Winsl{\o}w, Carl Ultraproducts, QWEP von Neumann algebras, and the Effros-Marechal topology. J. Reine Angew. Math. 715 (2016), 231-250.

\item Haagerup, Uffe Group $C^*$-algebras without the completely bounded approximation property. J. Lie Theory 26 (2016), no. 3, 861-887.

\item Haagerup, Uffe; de Laat, Tim Simple Lie groups without the approximation property II. Trans. Amer. Math. Soc. 368 (2016), no. 6, 3777-3809.

\item Haagerup, Uffe; Knudby, S{\o}ren The weak Haagerup property II: Examples. Int. Math. Res. Not. IMRN 2015, no. 16, 6941-6967.

\item Haagerup, Uffe; Musat, Magdalena An asymptotic property of factorizable completely positive maps and the Connes embedding problem. Comm. Math. Phys. 338 (2015), no. 2, 721-752.

\item Haagerup, S{\o}ren; Haagerup, Uffe; Ramirez-Solano, Maria A computational approach to the Thompson group F. Internat. J. Algebra Comput. 25 (2015), no. 3, 381-432.

\item Haagerup, Uffe; Knudby, S{\o}ren A Levy-Khinchin formula for free groups. Proc. Amer. Math. Soc. 143 (2015), no. 4, 1477-1489.

\item Haagerup, Uffe; Thorbj{\o}rnsen, Steen On the free gamma distributions. Indiana Univ. Math. J. 63 (2014), no. 4, 1159-1194.

\item Haagerup, Uffe Quasitraces on exact $C^*$-algebras are traces. C. R. Math. Acad. Sci. Soc. R. Can. 36 (2014), no. 2-3, 67-92.

\item Haagerup, U. Applications of random matrices to operator algebra theory. XVIIth International Congress on Mathematical Physics, 67, World Sci. Publ., Hackensack, NJ, 2014.

\item Ando, Hiroshi; Haagerup, Uffe Ultraproducts of von Neumann algebras. J. Funct. Anal. 266 (2014), no. 12, 6842-6913.

\item Haagerup, Uffe; Schlichtkrull, Henrik Inequalities for Jacobi polynomials. Ramanujan J. 33 (2014), no. 2, 227-246.

\item Haagerup, Uffe; M{\o}ller, S{\o}ren The law of large numbers for the free multiplicative convolution. Operator algebra and dynamics, 157-186, Springer Proc. Math. Stat., 58, Springer, Heidelberg, 2013.

\item Haagerup, Uffe; de Laat, Tim Simple Lie groups without the approximation property. Duke Math. J. 162 (2013), no. 5, 925-964.

\item Haagerup, Uffe; M{\o}ller, S{\o}ren Radial multipliers on reduced free products of operator algebras. J. Funct. Anal. 263 (2012), no. 8, 2507-2528.

\item Haagerup, Uffe; Thorbj{\o}rnsen, Steen Asymptotic expansions for the Gaussian unitary ensemble. Infin. Dimens. Anal. Quantum Probab. Relat. Top. 15 (2012), no. 1, 1250003, 41 pp.

\item Asaeda, Marta; Haagerup, Uffe Fusion rules on a parametrized series of graphs. Pacific J. Math. 253 (2011), no. 2, 257-288.

\item Haagerup, Uffe; Picioroaga, Gabriel New presentations of Thompson's groups and applications. J. Operator Theory 66 (2011), no. 1, 217-232.

\item Haagerup, Uffe; Musat, Magdalena Factorization and dilation problems for completely positive maps on von Neumann algebras. Comm. Math. Phys. 303 (2011), no. 2, 555-594.

\item Haagerup, U.; Steenstrup, T.; Szwarc, R. Schur multipliers and spherical functions on homogeneous trees. Internat. J. Math. 21 (2010), no. 10, 1337-1382.

\item Haagerup, Uffe; Kemp, Todd; Speicher, Roland Resolvents of R-diagonal operators. Trans. Amer. Math. Soc. 362 (2010), no. 11, 6029-6064.

\item Haagerup, Uffe; Junge, Marius; Xu, Quanhua A reduction method for noncommutative $L_p$-spaces and applications. Trans. Amer. Math. Soc. 362 (2010), no. 4, 2125-2165.

\item Haagerup, Uffe; Musat, Magdalena Classification of hyperfinite factors up to completely bounded isomorphism of their preduals. J. Reine Angew. Math. 630 (2009), 141-176.

\item Haagerup, Uffe; Schultz, Hanne Invariant subspaces for operators in a general II1-factor. Publ. Math. Inst. Hautes Etudes Sci. No. 109 (2009), 19-111.

\item Haagerup, Uffe; Musat, Magdalena The Effros-Ruan conjecture for bilinear forms on $C^*$-algebras. Invent. Math. 174 (2008), no. 1, 139-163.

\item Haagerup, Uffe; Musat, Magdalena On the best constants in noncommutative Khintchine-type inequalities. J. Funct. Anal. 250 (2007), no. 2, 588-624.

\item Haagerup, Uffe; Schultz, Hanne Brown measures of unbounded operators affiliated with a finite von Neumann algebra. Math. Scand. 100 (2007), no. 2, 209-263.

\item Haagerup, Uffe; Kadison, Richard V.; Pedersen, Gert K. Means of unitary operators, revisited. Math. Scand. 100 (2007), no. 2, 193-197.

\item Haagerup, Uffe; Schultz, Hanne; Thorbj{\o}rnsen, Steen A random matrix approach to the lack of projections in $C^*$red(F2). Adv. Math. 204 (2006), no. 1, 1-83.

\item Haagerup, Uffe; Thorbj{\o}rnsen, Steen A new application of random matrices: Ext($C^*$red(F2)) is not a group. Ann. of Math. (2) 162 (2005), no. 2, 711-775.

\item Aagaard, Lars; Haagerup, Uffe Moment formulas for the quasi-nilpotent DT-operator. Internat. J. Math. 15 (2004), no. 6, 581-628.

\item Dykema, Ken; Haagerup, Uffe Invariant subspaces of the quasinilpotent DT-operator. J. Funct. Anal. 209 (2004), no. 2, 332-366.

\item Dykema, Ken; Haagerup, Uffe DT-operators and decomposability of Voicu \newline -lescu's circular operator. Amer. J. Math. 126 (2004), no. 1, 121-189.

\item Haagerup, Uffe; Thorbj{\o}rnsen, Steen Random matrices with complex Gaussian entries. Expo. Math. 21 (2003), no. 4, 293-337.

\item Haagerup, U.; Rosenthal, H. P.; Sukochev, F. A. Banach embedding properties of non-commutative $L_p$-spaces. Mem. Amer. Math. Soc. 163 (2003), no. 776, vi+68 pp.

\item Haagerup, U. Random matrices, free probability and the invariant subspace problem relative to a von Neumann algebra. Proceedings of the International Congress of Mathematicians, Vol. I (Beijing, 2002), 273-290, Higher Ed. Press, Beijing, 2002. 46-02 pp.

\item Dykema, K.; Haagerup, U. Invariant subspaces of Voiculescu's circular operator. Geom. Funct. Anal. 11 (2001), no. 4, 693-741.

\item Haagerup, Uffe; Rosenthal, Haskell P.; Sukochev, Fedor A. On the
  Banach-isomorphic classification of $L_p$ spaces of hyperfinite von
  Neumann algebras. C. R. Acad. Sci. Paris Ser. I Math. 331 (2000),
  no. 9, 691-695.

\item Haagerup, U.; Thorbj{\o}rnsen, S. Random matrices and non-exact
  $C^*$-algebras. $C^*$-algebras (M{\"u}nster, 1999), 71-91, Springer,
  Berlin, 2000.

\item Haagerup, Uffe; Larsen, Flemming Brown's spectral distribution
  measure for R-diagonal elements in finite von Neumann
  algebras. J. Funct. Anal. 176 (2000), no. 2, 331-367.

\item Haagerup, Uffe; Winsl{\o}w, Carl The Effros-Marechal topology in
  the space of von Neumann algebras. II. J. Funct. Anal. 171 (2000),
  no. 2, 401-431.

\item Haagerup, U.; Thorbj{\o}rnsen, S. Random matrices and K-theory
  for exact $C^*$-algebras. Doc. Math. 4 (1999), 341-450.

\item Asaeda, M.; Haagerup, U. Exotic subfactors of finite depth with Jones indices $(5+\sqrt{13})/2$ and
 $(5+\sqrt{17})/2$. Comm. Math. Phys. 202 (1999), no. 1, 1-63.

\item Haagerup, Uffe; St{\o}rmer, Erling On maximality of entropy in
  finite von Neumann algebras. Operator algebras and operator theory
  (Shanghai, 1997), 99-109, Contemp. Math., 228, Amer. Math. Soc.,
  Providence, RI, 1998.

\item Haagerup, Uffe; Laustsen, Niels J. Weak amenability of
  $C^*$-algebras and a theorem of Goldstein. Banach algebras '97
  (Blaubeuren), 223-243, de Gruyter, Berlin, 1998.

\item Dykema, Ken; Haagerup, Uffe; R{\o}rdam, Mikael Correction to:
  "The stable rank of some free product $C^*$-algebras''. Duke
  Math. J. 94 (1998), no. 1, 213.

\item Haagerup, Uffe; Winsl{\o}w, Carl The Effros-Marechal topology in
  the space of von Neumann algebras. Amer. J. Math. 120 (1998), no. 3,
  567-617.

\item Haagerup, Uffe; St{\o}rmer, Erling Maximality of entropy in
  finite von Neumann algebras. Invent. Math. 132 (1998), no. 2,
  433-455.

\item Haagerup, Uffe Orthogonal maximal abelian $*$-subalgebras of the
  $n\times n$ matrices and cyclic n-roots. Operator algebras and
  quantum field theory (Rome, 1996), 296-322, Int. Press, Cambridge,
  MA, 1997.

\item Dykema, Ken; Haagerup, Uffe; R{\o}rdam, Mikael The stable rank
  of some free product $C^*$-algebras. Duke Math. J. 90 (1997), no. 1,
  95-121.

\item Haagerup, Uffe On Voiculescu's R- and S-transforms for free
  non-commuting random variables. Free probability theory (Waterloo,
  ON, 1995), 127-148, Fields Inst. Commun., 12, Amer. Math. Soc.,
  Providence, RI, 1997.

\item Bisch, Dietmar; Haagerup, Uffe Composition of subfactors: new
  examples of infinite depth subfactors. Ann. Sci. Ecole
  Norm. Sup. (4) 29 (1996), no. 3, 329-383.

\item Haagerup, Uffe; St{\o}rmer, Erling Positive projections of von Neumann algebras onto JW-algebras. Proceedings of the XXVII Symposium on Mathematical Physics (Torun, 1994). Rep. Math. Phys. 36 (1995), no. 2-3, 317-330.

\item Haagerup, Uffe; Itoh, Takashi Grothendieck type norms for
  bilinear forms on $C^*$-algebras. J. Operator Theory 34 (1995),
  no. 2, 263-283.

\item Haagerup, Uffe; R{\o}rdam, Mikael Perturbations of the rotation
  $C^*$-algebras and of the Heisenberg commutation relation. Duke
  Math. J. 77 (1995), no. 3, 627-656.

\item Haagerup, Uffe; St{\o}rmer, Erling Subfactors of a factor of
  type III? which contain a maximal centralizer. Internat. J. Math. 6
  (1995), no. 2, 273-277.

\item Haagerup, Uffe Principal graphs of subfactors in the index range
  $4<[M:N]<3+\sqrt{2}$. Subfactors (Kyuzeso, 1993), 1-38, World
  Sci. Publ., River Edge, NJ, 1994.

\item Haagerup, Uffe; Zsido, Laszlo Resolvent estimate for Hermitian
  operators and a related minimal extrapolation problem. Acta
  Sci. Math. (Szeged) 59 (1994), no. 3-4, 503-524.

\item Haagerup, Uffe; St{\o}rmer, Erling Pointwise inner automorphisms
  of injective factors. J. Funct. Anal. 122 (1994), no. 2, 307-314.

\item Haagerup, Uffe; Kraus, Jon Approximation properties for group
  $C^*$-algebras and group von Neumann
  algebras. Trans. Amer. Math. Soc. 344 (1994), no. 2, 667-699.

\item Haagerup, Uffe; R{\o}rdam, Mikael $C^*$-algebras of unitary rank
  two. J. Operator Theory 30 (1993), no. 1, 161-171.

\item Haagerup, Uffe; Pisier, Gilles Bounded linear operators between
  $C^*$-algebras. Duke Math. J. 71 (1993), no. 3, 889-925.

\item Haagerup, Uffe; de la Harpe, Pierre The numerical radius of a
  nilpotent operator on a Hilbert space. Proc. Amer. Math. Soc. 115
  (1992), no. 2, 371-379.

\item Anderson, Joel; Blackadar, Bruce; Haagerup, Uffe Minimal
  projections in the reduced group $C^*$-algebra of Zn*Zm. J. Operator
  Theory 26 (1991), no. 1, 3-23.

\item Haagerup, Uffe On convex combinations of unitary operators in
  $C^*$-algebras. Mappings of operator algebras (Philadelphia, PA,
  1988), 1-13, Progr. Math., 84, Birkhauser Boston, Boston, MA, 1990.

\item Haagerup, Uffe; St{\o}rmer, Erling Automorphisms which preserve
  unitary equivalence classes of normal states. Operator theory:
  operator algebras and applications, Part 1 (Durham, NH, 1988),
  531-537, Proc. Sympos. Pure Math., 51, Part 1, Amer. Math. Soc.,
  Providence, RI, 1990.

\item Haagerup, Uffe; St{\o}rmer, Erling Equivalence of normal states
  on von Neumann algebras and the flow of weights. Adv. Math. 83
  (1990), no. 2, 180-262.

\item Haagerup, Uffe; St{\o}rmer, Erling Pointwise inner automorphisms
  of von Neumann algebras. With an appendix by Colin
  Sutherland. J. Funct. Anal. 92 (1990), no. 1, 177-201.

\item Haagerup, Uffe; Pisier, Gilles Factorization of analytic functions with values in noncommutative L1-spaces and applications. Canad. J. Math. 41 (1989), no. 5, 882-906.

\item Cowling, Michael; Haagerup, Uffe Completely bounded multipliers of the Fourier algebra of a simple Lie group of real rank one. Invent. Math. 96 (1989), no. 3, 507-549.

\item Haagerup, Uffe The injective factors of type ${\rm III}_\lambda,\;0<\lambda<1$. Pacific J. Math. 137 (1989), no. 2, 265-310.

\item Cowling, M.; Haagerup, U.; Howe, R. Almost L2 matrix coefficients. J. Reine Angew. Math. 387 (1988), 97-110.

\item Haagerup, Uffe A new upper bound for the complex Grothendieck constant. Israel J. Math. 60 (1987), no. 2, 199-224.

\item Haagerup, Uffe Connes' bicentralizer problem and uniqueness of the injective factor of type III1. Acta Math. 158 (1987), no. 1-2, 95-148.

\item Haagerup, Uffe . Proceedings of the nineteenth Nordic congress of mathematicians (Reykjav\'{\i}k, 1984), 60-77, V\'{\i}sindaf\'el. \'{I}sl., XLIV, Icel. Math. Soc., Reykjav\'{\i}k, 1985.

\item Haagerup, Uffe Injectivity and decomposition of completely
  bounded maps. Operator algebras and their connections with topology
  and ergodic theory (Busteni, 1983), 170-222, Lecture Notes in Math.,
  1132, Springer, Berlin, 1985.

\item Connes, Alain; Haagerup, Uffe; St{\o}rmer, Erling Diameters of
  state spaces of type III factors. Operator algebras and their
  connections with topology and ergodic theory (Busteni, 1983),
  91-116, Lecture Notes in Math., 1132, Springer, Berlin, 1985.

\item Haagerup, Uffe A new proof of the equivalence of injectivity and hyperfiniteness for factors on a separable Hilbert space. J. Funct. Anal. 62 (1985), no. 2, 160-201.

\item Effros, Edward G.; Haagerup, Uffe Lifting problems and local reflexivity for $C^*$-algebras. Duke Math. J. 52 (1985), no. 1, 103-128.

\item Haagerup, Uffe The Grothendieck inequality for bilinear forms on $C^*$-algebras. Adv. in Math. 56 (1985), no. 2, 93-116.

\item De Canniere, Jean; Haagerup, Uffe Multipliers of the Fourier algebras of some simple Lie groups and their discrete subgroups. Amer. J. Math. 107 (1985), no. 2, 455-500.

\item Haagerup, Uffe; Hanche-Olsen, Harald Tomita--Takesaki theory for
  Jordan algebras. J. Operator Theory 11 (1984), no. 2, 343-364.

\item Haagerup, Uffe; Zsido, Laszlo Sur la propriete de Dixmier pour
  les $C^*$-algebres. (French) [On the Dixmier property for
  $C^*$-algebras] C. R. Acad. Sci. Paris Ser. I Math. 298 (1984),
  no. 8, 173-176.

\item Haagerup, U. All nuclear $C^*$-algebras are
  amenable. Invent. Math. 74 (1983), no. 2, 305-319.

\item Haagerup, Uffe Solution of the similarity problem for cyclic
  representations of $C^*$-algebras. Ann. of Math. (2) 118 (1983),
  no. 2, 215-240.

\item Haagerup, Uffe The best constants in the Khintchine
  inequality. Studia Math. 70 (1981), no. 3, 231-283 (1982).

\item Haagerup, Uffe The reduced $C^*$-algebra of the free group on
  two generators. 18th Scandinavian Congress of Mathematicians
  (Aarhus, 1980), pp. 321-335, Progr. Math., 11, Birkheuser, Boston,
  Mass., 1981.

\item Haagerup, Uffe; Skau, Christian F. Geometric aspects of the
  Tomita-Takesaki theory. II. Math. Scand. 48 (1981), no. 2, 241-252.

\item Haagerup, Uffe; Munkholm, Hans J. Acta Math. 147 (1981),
  no. 1-2, 1-11.

\item Haagerup, Uffe $L_p$-spaces associated with an arbitrary von
  Neumann algebra. Algebres d'operateurs et leurs applications en
  physique mathematique (Proc. Colloq., Marseille, 1977), pp. 175-184,
  Colloq. Internat. CNRS, 274, CNRS, Paris, 1979.

\item Haagerup, Uffe Operator-valued weights in von Neumann
  algebras. II. J. Funct. Anal. 33 (1979), no. 3, 339-361.

\item Haagerup, Uffe A density theorem for left Hilbert
  algebras. Algebres d'operateurs (Sem., Les Plans-sur-Bex, 1978),
  pp. 170-179, Lecture Notes in Math., 725, Springer, Berlin, 1979.

\item Haagerup, Uffe Operator-valued weights in von Neumann
  algebras. I. J. Funct. Anal. 32 (1979), no. 2, 175-206.

\item Haagerup, Uffe The best constants in the Khintchine
  inequality. Proceedings of the International Conference on Operator
  Algebras, Ideals, and their Applications in Theoretical Physics
  (Leipzig, 1977), pp. 69-79, Teubner, Leipzig, 1978.

\item Haagerup, Uffe On the dual weights for crossed products of von Neumann algebras. II. Application of operator-valued weights. Math. Scand. 43 (1978/79), no. 1, 119-140.

\item Haagerup, Uffe On the dual weights for crossed products of von Neumann algebras. I. Removing separability conditions. Math. Scand. 43 (1978/79), no. 1, 99-118.

\item Haagerup, Uffe An example of a nonnuclear $C^*$-algebra, which has the metric approximation property. Invent. Math. 50 (1978/79), no. 3, 279-293.

\item Haagerup, Uffe Les meilleures constantes de l'inegalite de
  Khintchine. (French) C. R. Acad. Sci. Paris Ser. A-B 286 (1978),
  no. 5, A259-A262.

\item Bratteli, Ola; Haagerup, Uffe Unbounded derivations and
  invariant states. Comm. Math. Phys. 59 (1978), no. 1, 79-95.

\item Haagerup, Uffe An example of a weight with type III
  centralizer. Proc. Amer. Math. Soc. 62 (1977), no. 2, 278-280.

\item Haagerup, Uffe Operator valued weights and crossed
  products. Symposia Mathematica, Vol. XX (Convegno sulle Algebre
  $C^*$ e loro Applicazioni in Fisica Teorica, Convegno sulla Teoria
  degli Operatori Indice e Teoria K, INDAM, Roma, 1974),
  pp. 241-251. Academic Press, London, 1976.

\item Haagerup, Uffe The standard form of von Neumann
  algebras. Math. Scand. 37 (1975), no. 2, 271-283.

\item Haagerup, Uffe Normal weights on
  $W^*$-algebras. J. Funct. Anal. 19 (1975), 302-317.

\end{itemize}

\textbf{Acknowledgements.} The authors would like to sincerely thank
the sons of Uffe, Peter and S{\o}ren, for providing us several
memorable photos of Uffe and for their valuable suggestions improving
the biography of Uffe. The authors also wish to express their gratitude to
Professor Alain Connes for supplying the photograph shown in Figure
7.

\bibliographystyle{amsplain}

\begin{thebibliography}{99}

\bibitem{HA} M. Asaeda and U. Haagerup, \textit{Exotic subfactors of
 finite depth with Jones indices $(5+\sqrt{13})/2$ and
 $(5+\sqrt{17})/2$}, Comm. Math. Phys. \textbf{202} (1999), no. 1,
 1--63.

\bibitem{LB_MB} A. Le Boudec and N. Matte Bon, \textit{Subgroup
 dynamics and $C^*$-simplicity of groups of homeomorphisms},
 ArXiv:1605.01651 (2016)

\bibitem{COW} P.-A. Cherix, M. Cowling, P. Jolissaint, P. Julg, and
  A. Valette, \textit{Groups with the Haagerup property. Gromov's
    a-T-menability}, Progress in Mathematics, 197.
Birkh\"auser Verlag, Basel, 2001.

\bibitem{ALAIN3} A. Connes, \textit{Une classification des facteurs de type ${\rm III}$. (French)}, Ann. Sci. \'Ecole Norm. Sup. (4) \textbf{6} (1973), 133--252.

\bibitem{Co} A. Connes, \textit{Factors of type III$_1$, property
 $L'_{\lambda}$ and closure of inner automorphisms},
J. Operator Theory \textbf{14} (1985), 189--211.

\bibitem{ALAIN} A. Connes, \textit{Uffe Haagerup},
  http://noncommutativegeometry.blogspot.co.uk/2015/07/uffe-haagerup.html,
  Saturday, July 18, 2015.

\bibitem{ALAIN2} A. Connes, V. Jones, M. Musat, M. R{\o}rdam,
  \textit{Uffe Haagerup in memoriam}, Notices
  Amer. Math. Soc. \textbf{63} (2016), no. 1, 48--49.

\bibitem{HHR} S. Haagerup, U. Haagerup, and M. Ramirez-Solano,
  \textit{Computational explorations of the
 Thompson group T for the amenability problem of F}.
 ArXiv:1705.00198.

\bibitem{first} U. Haagerup, \textit{The standard form of von Neumann
    algebras}, Math. Scand. \textbf{37} (1975), no. 2, 271--283.

\bibitem{second}  U. Haagerup, \textit{Normal weights on
    $W^*$-algebras}, J. Funct. Anal. 19 (1975), 302-317.

\bibitem{msm1} U. Haagerup, \textit{An example of a nonnuclear
    $C^*$-algebra, which has the metric approximation property},
  Invent. Math. 50 (1978/79), no. 3, 279--293.

\bibitem{haa2} U. Haagerup, \textit{Connes' bicentralizer problem
 and uniqueness of the injective factor of type III$_1$},
Acta Math. \textbf{158} (1987), no. 1-2, 95--148.

\bibitem{msm2} U. Haagerup, \textit{Principal graphs of subfactors in
    the index range $4<[M:N]<3+\sqrt{2}$}, Subfactors (Kyuzeso, 1993),
  1-38, World Sci. Publ., River Edge, NJ, 1994.

\bibitem{haa3} U. Haagerup, \textit{On Voiculscu's $R$- and
 $S$-transforms for free non-commuting random variables},
Free probability theory (Waterloo, ON, 1995), 127-148, Fields
Inst. Commun. \textbf{12}, Amer. Math. Soc., Providence, RI, 1997.

\bibitem{H_KO} U. Haagerup and K. Knudsen Olesen, \textit{Non-inner
 amenability of the Thompson groups $T$ and $V$}, Journ. Funct.
 Anal. \textbf{272} (2017), 4838--4852.

\bibitem{H_dL1} U. Haagerup and T. de~Laat, \textit{Simple Lie
 Groups without the Approximation Property}, Duke
Math. J. \textbf{162} (2013), 925--964.

\bibitem{H_dL2} U. Haagerup and T. de~Laat, \textit{Simple Lie Groups
    without the Approximation Property II}, Trans. Amer.
 Math. Soc. \textbf{368} (2016), 3777--3809.

\bibitem{HaMu} U. Haagerup and M. Musat, \textit{The Effros-Ruan
    conjecture for bilinear forms on $C^*$-algebras},
  Invent. Math. \textbf{174} (2008), no. 1, 139-163.

\bibitem{HS}
U. Haagerup and H. Schultz, \textit{Invariant subspaces for operators
  in a general ${\rm II}_1$-factor},
Publ. Math. Inst. Hautes \'Etudes Sci. No. \textbf{109} (2009), 19--111.

\bibitem{ht1}
U. Haagerup and S. Thorbj{\o}rnsen, {\sl Random Matrices with Complex
  Gaussian Entries}, Expositiones Math. \textbf{21} (2003), 293--337.

\bibitem{ht5}
U. Haagerup and S. Thorbj{\o}rnsen, \textit{A new application of random
  matrices: ${\rm Ext}(C^*_{\rm red}(F_2))$ is not a group}, Annals of
Math. \textbf{162} (2005), 711--775.

\bibitem{interview} J. Hjelmborg, \textit{Interview med Uffe Haagerup},
  Matilde (newletter of the Danish Mathematical Society) \textbf{12}
  (2002).

\bibitem{L_dlS} V. Lafforgue and M. de~la~Salle, \textit{
 Noncommutative $L^p$-spaces without the completely bounded
 approximation property}, Duke Math. J. \textbf{160} (2011), 71--116.

\bibitem{NiSp} A. Nica and R. Speicer, \textit{Lectures on the
    Combinatorics of Free Probability}, London Math.\ Soc.\ Lecture
  Note Series \textbf{335}, Cambridge University Press (2006).

\bibitem{Pi_Sh} G. Pisier and D. Shlyakhtenko, \textit{
 Grothendieck's theorem for operator spaces},
Invent. Math. \textbf{150} (2002), 185--217.

\bibitem{sigmund} P. Sigmund and U. Haagerup, \textit{Bethe stopping
    theory for a harmonic oscillator and Bohr's oscillator model of
    atomic stopping}, Physical Review A (General Physics), \textbf{34}
  (1986), Issue 2, 892--910.

\bibitem{Vo} D. Voiculescu, \textit{Addition of certain
    non-commutative random variables}, J. Funct. Anal.\ \textbf{66}
  (1986), 323-346.

\bibitem{Wi} E. P. Wigner, \textit{On the distribution of the roots of
    certain random matrices}, Ann. Math.\ \textbf{67} (1958), 325-327.

\bibitem{cv} \url{http://www.math.ku.dk/~haagerup/index.php?show=cv}

\bibitem{mr} \url{http://www.ams.org/mathscinet/}

\bibitem{zbl} \url{https://zbmath.org/}

\bibitem{wos} \url{https://www.webofknowledge.com/}

\bibitem{scopus} \url{https://www.scopus.com/}

\bibitem{wik} Wikipedia: \url{https://en.wikipedia.org/wiki/Uffe_Haagerup}


\end{thebibliography}

\end{document}